\documentclass{amsart}
\usepackage{stmaryrd}
\usepackage{amssymb}
\usepackage{amsmath}
\usepackage{amsthm}
\usepackage{amsbsy}
\usepackage[dvips]{graphics}
\usepackage{amscd}
\usepackage{bm}

\theoremstyle{plain}
\newtheorem{theorem}{Theorem}[section]

\newtheorem{lemma}[theorem]{Lemma}
\newtheorem{corollary}[theorem]{Corollary}

\theoremstyle{definition}

\theoremstyle{remark}

\newtheorem*{acknowledgements}{Acknowledgements}

\newcommand{\A}{\mathfrak{A}}

\renewcommand{\H}{\mathfrak{H}}

\newcommand{\C}{{\mathbb{C}}}
\newcommand{\R}{{\mathbb{R}}}

\begin{document} 

\title[Scalar curvature and minimal hypersufaces]{Positive scalar
curvature and minimal hypersurfaces}

\author{Harish Seshadri}
\email{harish@isibang.ac.in}
\address{Stat-Math Unit,\\
Indian Statistical Institute,\\
Bangalore, India}

\subjclass{Primary 53C21}
\date{\today}

\begin{abstract}
We show that the minimal hypersurface method of Schoen and Yau can be
used for the ``quantitative'' study of positive scalar curvature. More
precisely, we show that if a manifold admits a metric $g$ with $s_g \ge
\vert T \vert$ or $s_g \ge \vert W \vert$, where $s_g$ is the scalar
curvature of of $g$, $T$ any 2-tensor on $M$ and $W$ the Weyl tensor
of $g$, then any closed orientable stable minimal (totally geodesic in
the second case) hypersurface also admits a metric with the corresponding
positivity of scalar curvature.  A corollary about the topology of such
hypersurfaces is proved in a special situation.
\end{abstract}

\maketitle
\section{Introduction}
In the study of the topology and geometry of Riemannian manifolds with
positive scalar curvature, two approaches have played a fundamental
role. The first one was initiated by Lichnerowicz in 1961 ~\cite{lic}
and is based on the Bochner-Weitzenbock formula for the Dirac operator
acting on spinors. This method was greatly developed in a series of
papers by Gromov and Lawson ~\cite{gro}. The second method made its
appearance in the paper of Schoen and Yau in 1979 ~\cite{sch} and used
the theory of minimal hypersurfaces. One of the applications of these
methods was to the question of which manifolds admit complete metrics
of positive scalar curvature.

It is natural to seek ``quantitative'' versions of the above
``qualitative'' question. Our notion of ``quantitative'' is the
following: let $(M,g)$ be a Riemannian manifold and $s$ (or $s_g$, if
we want to emphasize the dependence on $g$) the scalar curvature of
$g$. If $f: M \rightarrow {\Bbb R}$ is a function (depending on $g$
and not necessarily smooth), we are interested in the topological
implications of the condition $s \ge f$. In this paper we consider
functions $f=f_g$ that are of {\it conformal weight -2} in the
terminology of ~\cite{gur}. This means that $f$ is a function which
depends on the metric $g$ in the following way: $f_{u^2g}=u^{-2}f_g$.

It can be shown (cf. ~\cite{gur}) that, if $f$ is Lipschitz continuous,
the conformal class of every metric contains a $C^{2, \alpha}$ metric for
which $s>f$, $s<f$ or $s=f$ and these possibilities are mutually exclusive.
Hence the condition `` $s \ge f$'' is actually a conformal notion.
  
We consider three examples of such $f$:

1) $f= \vert T \vert_g$ where $T$ is 2-tensor on $M$ and $\vert T
\vert_g$ denotes the pointwise norm of $T$ with respect to $g$.

2) $f= \vert W_g \vert_g$ where $W_g$ is the Weyl tensor 
of $g$ and $\vert W \vert_g$ again denotes pointwise norm.

3) Let $i: \Sigma^{n-1} \rightarrow M^{n}$ be a hypersurface, with $i$
being the inclusion map. Fix a Riemannian metric $h$ on $\Sigma$ and
let $g$ be a fixed extension of $h$ to $M$. Let $u$ be any smooth
positive function on $\Sigma$ and let $v$ denote any extension of $u$
to a neighbourhood of $\Sigma$. Denote the second fundamental form and
mean curvature (with respect to a fixed normal vector field) of
$\Sigma$ in $(M,g)$ by $B_g$ and $H_g$ respectively.  If $L_{u^2 h} =
\Bigl ( B_{v^2g}-{\frac {H_{v^2g}}{n-1}}u^2h \Bigr )$ denotes the
trace-free part of the second fundamental form of $(\Sigma, u^2h)$ in
$(M,v^2g)$, then it can be checked that $L$ does not depend on the
extension $v$. Also,
$$ \vert L_{u^2h } \vert^2 _{u^2 h} = u^{-2}\vert L_{h} \vert^2
_{h}.$$ Hence, for $h$ and $g$ as above, the function $f= \vert L_h
\vert^2_h$ is of conformal weight $-2$. \\ 

When $T$ is a closed 2-form on $M$, Gursky and LeBrun ~\cite{gur} used
a variation of the spinor method to obtain obstructions to the
existence of metrics with $s \ge \vert T \vert $ and applied this for
the calculation of the Yamabe invariant of certain $4$-manifolds. The
purpose of this note is to show that the minimal hypersurface method
also yields information in the context of quantitative positive scalar
curvature and derive some geometric corollaries. More precisely, we have

\begin{theorem}\label{mai1}
Let $(M,g)$ be an orientable Riemannian $n$-manifold with $n \ge 3$.
Assume that $s_g \ge \vert T \vert _g$ where $T$ is a smooth 2-tensor
on $M$. If $ i: \Sigma \rightarrow M$ is a closed, orientable, stable
minimal hypersurface in $M$, then 
one of the two possibilities hold (where $L_h$ is interpreted as in Example
3)):
\begin{enumerate}
\item[{(i)}] $\Sigma$ admits a $C^{2, \alpha}$ metric $h$ conformal to
$i^*(g)$ with $s_h > \vert i^{*}(T) \vert _h + \vert L_h \vert^2_h$.
\item[{(ii)}] The induced
metric $h$ on $\Sigma$ satisfies $s_h \equiv \vert i^{*}(T) \vert _h +
\vert B_h \vert^2_h$.
\end{enumerate}
\end{theorem}

\noindent A similar result holds for $f= \vert W \vert$ :
\begin{theorem}\label{mai2}
Let $(M,g)$ be an orientable Riemannian $n$-manifold with $n \ge
5$. Assume that $s_g \ge c \vert W_g \vert _g$, where $c \in {\Bbb
R}^+$. If $ i: \Sigma \rightarrow M$ is a closed, orientable, stable
totally geodesic hypersurface in $M$, then one of the two
possibilities hold:
\begin{enumerate}
\item[{(i)}] $\Sigma$ admits a $C^{2, \alpha}$
metric $h$ conformal to $i^*(g)$ with $s_h > c \vert W_h \vert _h $.
\item[{(ii)}] The induced metric $h$ on $\Sigma$ satisfies
$s_h \equiv c \vert W_h \vert _h $.
\end{enumerate}
\end{theorem}
In (ii) of both these theorems, the proofs give more information than what has
been stated above. We refer the reader to the proofs for details.
 
Combining Theorem \ref{mai2} with the B\"ochner formula for harmonic
$2$-forms on $4$-manifolds and , we obtain the
following corollary.
\begin{corollary} \label{c1}
Let $(M,g)$ be an orientable Riemannian 5-manifold with $ s \ge 2
\sqrt 6 \vert W \vert$. If $i: \Sigma \rightarrow M$ is a closed
orientable stable totally geodesic hypersurface, then either $H^2(
\Sigma, {\Bbb R})=0$ or $(\Sigma, i^*(g))$ is a self-dual K\"ahler surface
of non-negative scalar curvature. 
\end{corollary}
Note that a theorem of Derdzinski ~\cite{der} states that any compact
self-dual K\"ahler surface is locally symmetric.

As examples of $5$-manifolds $(M,g)$ satisfying the hypotheses of
Corollary \ref{c1} we can consider $(M,g)=(\Sigma \times S^1, h+d
\theta^2)$, where $(\Sigma, h)$ is (upto scaling) either $S^4$ with
its standard metric, a closed flat 4-manifold or $ \C P^2$ with the
Fubini-Study metric. It is clear that the condition $ s \ge 2 \sqrt 6
\vert W \vert$ is satisfied for the cases of $S^4$ and flat
4-manifolds because the Weyl tensor of a Riemannian product of a
constant sectional curvature space with $S^1$ is zero. The
verification for the $ \C P^2$ case is the remark on Page 6. Note that
for this case $ s \equiv 2 \sqrt 6 \vert W \vert = 2 \sqrt 6 \vert W^+
\vert \neq 0$, where $W^+$ is the self-dual part of the Weyl tensor.
It is well-known that $s= 2 \sqrt 6 \vert W^+ \vert$ for any K\"ahler
surface $(\Sigma, h)$ with non-negative scalar curvature and with the
natural orientation, cf. ~\cite{der}. However, according to Corollary
\ref{c1}, only the self-dual ones can possibly arise as stable totally
geodesic hypersurfaces in 5-manifolds satisfying $ s \ge 2 \sqrt 6
\vert W \vert$.

For the next corollary we note that if a manifold $M$ has $Ric \ge 0$,
then the second variation formula for the area functional (see
(\ref{sec}) in Section 2) implies that any stable minimal hypersurface
in $M$ is totally geodesic. Combining this observation with the
existence of orientable stable minimal hypersurfaces representing
codimension-one homology classes in $M$ (when $M$ is orientable and
dim $M \le 7$) and Corollary \ref{c1}, we get
\begin{corollary}
Let $(M,g)$ be an orientable Riemannian 5-manifold with $Ric \ge 0$
and $ s \ge 2 \sqrt 6 \vert W \vert$. Then every homology
class $\alpha$ in $H_4(M, {\Bbb R})$ is represented by a 
hypersurface $\Sigma _\alpha $ with $H^2(\Sigma _ \alpha, {\Bbb R})=0$ or
by a self-dual K\"ahler surface of non-negative scalar curvature.
\end{corollary}

\begin{acknowledgements}
I would like to thank Kazuo Akutagawa and Siddartha Gadgil for several
interesting and helpful discussions.

I would also like to thank the referee for various important comments. 
\end{acknowledgements}

\section{Proofs}
{\bf Convention}: For the rest of the paper we follow the following
convention. For any tensor $S$ on $M$ and any submanifold $i: \Sigma
\rightarrow M$, 
$$ we \ denote \ the \ tensor \ i^*(S) \ on \ \Sigma \ by \ \tilde S.$$

\begin{proof}[Proof of Theorem \ref{mai1}] 

The proof is based on that of Schoen and Yau in ~\cite{sch}.

The object of our study will be the {\it modified scalar curvature}
$\sigma (g,f)$ corresponding to a function $f_g$ of conformal weight
$-2$ defined by $$ \sigma (g,f)= s_g -f_g.$$ 

We begin by noting that, under a conformal change of the metric $g
\rightarrow \overline g=u^{\frac{4}{n-2}} g$, the modified scalar
curvature transforms by 
\begin{equation}
\sigma \bigl (g,f ) \rightarrow \sigma
\bigl ( {\overline g}, f \bigr ) = u^{\frac{-4}{n-2}} \sigma \bigl (g,f
\bigr ) -4{\frac{n-1}{n-2}} u^{- {\frac{n-2}{n +2}}} \triangle u. \label{tra}
\end{equation}
This follows immediately from the corresponding transformation law for
the usual scalar curvature $s$ cf. ~\cite{bes}.

The stability of the minimal surface $\Sigma$ implies that for any
smooth function $\phi$ on $\Sigma$, 
\begin{equation}
 -\int _ {\Sigma} \bigl(Ric(N,N) + \vert B \vert ^2 \bigr) \phi^2 + 
\int _{\Sigma} \vert \nabla \phi \vert ^2 \ge 0, \label{sec}
\end{equation}
where $N$ is a unit normal vector field on $\Sigma$, $Ric$ is the
Ricci tensor of $g$ and $B$ is the second fundamental form of $\Sigma$.

As in ~\cite{sch}, the Gauss curvature equations and the minimality
of $\Sigma$ can be used to write the above inequality can as
\begin{equation}\label{sec1} 
\int _{\Sigma} {\frac{s \phi^2}{2}} -\int _{\Sigma}
{\frac{s_{\tilde g} \phi^2}{2}} + {\frac{1}{2}} \int _{\Sigma} \vert B
\vert^2 \phi^2 \le \int_{\Sigma} \vert \nabla \phi \vert ^2, 
\end{equation}
where $s_{\tilde g}$ is the scalar curvature of $\Sigma$ with the induced
metric.  

To adapt (\ref{sec1}) for modified scalar curvature, we need
the following simple lemma.
\begin{lemma}\label{nor}
Let $S \in \otimes^r T^*M^n$ be a tensor on $M$ and $i: \Sigma^k
\rightarrow M$ a submanifold, with $i$ being the inclusion map. Then
$\vert \tilde S \vert _{\tilde g} \le \vert S \vert_g$ at any point of
$\Sigma$. 
\end{lemma}
\begin{proof}
Let $p \in \Sigma$. Let $\{ e_1,..,e_k \}$ be an orthonormal basis,
with respect to $\tilde g$, of $T_p \Sigma$, where $k=dim \ \Sigma$. 
Extend it to an orthonormal basis $\{ e_1,..,e_n \}$ of
$T_pM$ (we identify $i_*(e_i)$ with $e_i$ for $i=1,..,k$). We have
\begin{align} 
 \vert \tilde S \vert^2 _{\tilde g} &= \Sigma_k
\tilde S(e_{i_1},..,e_{i_r})^2  \notag \\
       &=  \Sigma_k S(i_*(e_{i_1}),..,i_*(e_{i_r}))^2, \notag \\
       &= \Sigma_k S(e_{i_1},..,e_{i_r})^2 \notag 
\end{align}
where $\Sigma_k$ denotes that the sum is
over all $(i_1,..,i_r)$ with $1 \le i_j \le k, \ j=1,..,r$. Hence we
get
$$ \vert \tilde S \vert^2 _{\tilde g} \le \Sigma_n
S(e_{i_1},..,e_{i_r})^2, $$ where the sum is now over all
$(i_1,..,i_r)$ with $1 \le i_j \le n, \ j=1,..,r$. But the right-hand
side of the above inequality is $\vert S \vert^2 _g$.
\end{proof}

We continue with the proof of the Theorem. Note the hypothesis on $s$
and the above Lemma imply that
\begin{equation}\label{mn} 
s \ge \vert T \vert \ge \vert \tilde T \vert_{\tilde g}.
\end{equation}
Combining this with the stability inequality (\ref{sec1}), we have

\begin{equation}\label{mo}
\int _{\Sigma} {\frac{ \vert \tilde T \vert_{\tilde g} \phi^2}{2}}
-\int _{\Sigma} {\frac{s_{\tilde g} \phi^2}{2}} + {\frac{1}{2}} \int
_{\Sigma} \vert B \vert^2 \phi^2 \le \int_{\Sigma} \vert \nabla \phi
\vert^2. 
\end{equation} 

Now consider the conformal function $f$ of weight $-2$ on $ \Sigma $
defined by $$f:= \vert \tilde T \vert + \vert L \vert ^2,$$ where 
$L$ is the trace-free part of the second fundamental form $B$ and we
interpret the conformal invariance of $\vert L \vert^2$ as in Example 3)
on Page 2. The associated modified scalar curvature is $ \sigma \bigl
(h, f \bigr )=s_h-f_h$.

Since $\Sigma$ is minimal, we have $L_{\tilde g}=B_{\tilde g}$ and we
can hence rewrite (\ref{mo}) as
\begin{equation}\label{mo2}
- {\frac {1}{2}} \int _{\Sigma} \sigma \bigl
(\tilde g, f \bigr ) \phi ^2  \le \int_{\Sigma} \vert \nabla \phi \vert^2.
\end{equation}
Let $ \triangle $ be the Laplace operator of $\tilde g$ and consider
the equation
\begin{equation}\label{aux}
 \triangle \phi = {\frac {(n-3)}{4(n-2)}} \sigma \bigl (\tilde g, f
\bigr ) \phi+ \lambda \phi.
\end{equation}
If this equation has a nontrivial solution $\phi$, we claim that
$\lambda \le 0$. Suppose $ \lambda >0$. 

As in ~\cite{sch},
we multiply the above equation by $\phi$,
integrate and use $(\ref{mo2})$ to get
\begin{equation}\label{deu}
 2{\frac {(n-2)}{(n-3)}} \int_\Sigma \vert \nabla \phi \vert ^2 =
-{\frac {1}{2}} \int _\Sigma \sigma \bigl (\tilde g, f \bigr ) 
\phi^2 - 2 \lambda {\frac {(n-2)}{(n-3)}} \int _\Sigma \phi ^2 <
\int _\Sigma \vert \nabla \phi \vert ^2 ,
\end{equation}
which is impossible since $\phi \not \equiv 0$.

Now define the operator
$$L=\triangle -{\frac {(n-3)}{4(n-2)}} \sigma \bigl ( \tilde g, f \bigr )$$ and let
$$\lambda= inf _{{\psi \in W^{1,2}} \atop {\Vert \psi \Vert _2=1}}<L \psi,
\psi>$$
and $u >0$ be the corresponding eigenfunction. Note that since
$\sigma \bigl (\tilde g, f \bigr )$ is, in general, Lipschitz continuous
but not smooth (at the zero locus of $\vert T \vert$), the best regularity we
can obtain for $u$ is that $u \in C^{2, \alpha}$ for any $0< \alpha <1$. This
is sufficient for our purposes. By definition, $u$ satisfies 
\begin{equation}
 L(u)= \lambda u \label{pos}
\end{equation} 
where $\lambda \le 0$. Multiplying the metric of
$\Sigma$ by $u^{\frac{4}{n-3}}$ and using (\ref{tra}), we see that the
modified scalar curvature of $\Sigma$ changes to
$$ -4{\frac {(n-2)}{(n-3)}} u^{-{\frac {(n-3)}{4(n-2)}}-1} L(u)$$ 
and hence is positive unless $\lambda =0$, by (\ref{pos}).
In this case, which is (i) of Theorem \ref{mai1}, the metric $h$ is given 
by $h=u^{\frac{4}{n-3}} \tilde g$.
 
If $\lambda=0$, then since $u \not \equiv 0$ satisfies (\ref{aux}) and
hence (\ref{deu}), we get $u \equiv$ constant and $ \sigma \bigl 
(\tilde g, f \bigr ) \equiv 0$, i.e. $s_{\tilde g}- 
\vert \tilde T \vert_{\tilde g} -\vert B \vert^2_{\tilde g} \equiv 0$. 
This completes the proof.

Here we indicate the other conclusions one can draw: note that 
$$s_g-s_{\tilde g}+\vert B \vert^2_{\tilde g} \ge \vert T \vert_g -
\vert \tilde T \vert_{\tilde g} \ge 0.$$ Going back to (\ref{sec1})
and taking $\phi \equiv$ constant, we see that $s_g \equiv s_{\tilde
g}- \vert B \vert^2_{\tilde g}$ and hence
$\vert T \vert_g \equiv \vert \tilde
T \vert_{\tilde g}$ by the above inequality.
\end{proof}

Before proceeding with the proof of Theorem \ref{mai2}, we point out
an issue involving the norm $\vert W \vert$. We look at $W$ both as an
element of $ End(\Lambda^2 T_p^*M)$ (for instance in Corollary \ref{c1})
and of $ \otimes^4 T_p^*M$ (in the claim below). These spaces have
norms which differ by a fixed scale. Since this scale is independent
of dimension, it is readily seen that, for the results that we prove,
this does not pose a problem.

\begin{proof}[Proof of Theorem \ref{mai2}]
The proof is very similar to the proof of Theorem \ref{mai1}. The only
extra point to note is that in order to get the inequality
corresponding to the second part of (\ref{mn}), we need to know that
if $\Sigma$ is totally geodesic, then $\vert W_{\tilde g}
\vert_{\tilde g} \le \vert W_g \vert_g $. This follows from the
following claim, since Lemma ~\ref{nor} implies that 
\begin{equation}\label{che}
\vert \tilde
{W_g} \vert_{\tilde g} \le \vert W_g \vert_g .
\end{equation}

CLAIM: $\vert W_{\tilde g} \vert_{\tilde g} \le \vert \tilde {W_g}
\vert_{\tilde g} $ with equality if and only if $W_{\tilde g}=\tilde W_g$

Fix a point $p \in \Sigma$ and let $R_{\tilde g}$ be the Riemann
curvature tensor of $(\Sigma, \tilde g)$ at $p$, regarded as an
element of $S^2 \Lambda ^2 T^*_p \Sigma$. $R_{\tilde g}$ decomposes
orthogonally (see, for instance, ~\cite{bes}) as
\begin{equation}\label{dec0}
R_{\tilde g}= {\frac {s_{\tilde g}}{2(n-1)(n-2)}}{\tilde g} \owedge \tilde g 
+ {\frac {1}{n-3}} z_{\tilde g} \owedge \tilde g +W_{\tilde g}, 
\end{equation}
where $z_{\tilde g}=Ric_{\tilde g} -{\frac {s_g}{n-1}} \tilde g$ is
the trace-free part of the Ricci tensor and $\owedge$ denotes the
Kulkarni-Nomizu product.  

We note that, if $\A= \{ h \owedge \tilde g : h \in S^2T^*_p \Sigma \}$,
then $ < W_{\tilde g}, u>=0$ for any $u \in \A$.

For convenience, write (\ref{dec0}) as
\begin{equation}\label{dec}
R_{\tilde g}= S_{\tilde g} + W_{\tilde g},
\end{equation}
where $S_{\tilde g}$ is the sum of the first two terms on the right side of
(\ref{dec0}).

On the other hand, since $\Sigma$ is totally geodesic, 
\begin{equation}
R_{\tilde g}= \tilde R_g= \tilde S_g + \tilde W_g, \label{dec2}
\end{equation}
where $R_g$, $S_g$ and $W_g$ are the corresponding
quantities for $g$.
Combining (\ref{dec}) and (\ref{dec2}), we get (all the norms
are with respect to $\tilde g$)
\begin{equation}\label{cha}
\vert \tilde W_g \vert^2 =  \vert W_{\tilde g} + S_{\tilde g}- \tilde S_g 
\vert ^2  
                       = \vert W_{\tilde g} \vert ^2 + \vert S_{\tilde g} 
- \tilde S_g \vert ^2.                       
\end{equation}
The second equality follows from the fact that since $S_g=h \owedge g$
for some $h$, we have $\tilde S = \tilde h \owedge \tilde g \in \A$.

Hence $ \vert \tilde W_g \vert \ge \vert W_{\tilde g} \vert $.

In the case of equality, we have $S_{\tilde g}= \tilde
S_g$. (\ref{dec}) and (\ref{dec2}) then imply that $W_{\tilde g}=
\tilde W_g$. This completes the proof of the claim.

In the $ \lambda =0$ case, corresponding to the last paragraph of the
proof of Theorem \ref{mai1}, we get $s_{\tilde g} - \vert W_{\tilde g}
\vert_{\tilde g} \equiv 0$. This is case (ii) in \ref{mai2}.

\noindent Proceeding as before, we also get: $s_g=s_{\tilde g}$ and $\vert W_g
\vert_g = \vert  W_{\tilde g} \vert_{\tilde g}$.
(\ref{che}) then implies that
$\vert W_{\tilde g} \vert_{\tilde g}= \vert \tilde W_ g \vert_{\tilde
g}$. Hence $\tilde S_g= S_ {\tilde g}$, by the proof of the claim.
\end{proof}

{\bf Remark:} Let us consider the simple case of a Riemannian product
$(M,g)=(\Sigma \times S^1, \tilde g + d\theta^2)$ and prove that the
equality $\vert W_{\tilde g} \vert_{\tilde g}= \vert  W_ g
\vert_g$ holds if and only if $(\Sigma, \tilde g)$ is Einstein.
To see this, we check when equality holds in both the claim and
(\ref{che}). Equality in the claim holds if and only if
$\tilde S_g =S_{\tilde g}$. Hence,
\begin{equation}\label{sta}
{\frac {s_g}{2n(n-1)}} \tilde g \owedge \tilde g + {\frac {1}{n-2}}
\tilde z_g \owedge \tilde g
 \ = \ {\frac {s_{\tilde g}}{2(n-1)(n-2)}}{\tilde g} \owedge \tilde g 
+ {\frac {1}{n-3}} z_{\tilde g} \owedge \tilde g.
\end{equation}
First note that $\tilde Ric_g = Ric_{\tilde g}$ because
$R_g(X,N,N,Y)=0$ for any $X, \ Y$ in $T_p \Sigma$ and $N$ normal to
$\Sigma$. 

Hence 
$$\tilde z_g \ = \ \tilde Ric_g -{\frac {s_g}{n}} \tilde g \ = \
Ric_{\tilde g} -{\frac {s_g}{n-1}} \tilde g + {\frac {s_g}{n(n-1)}}
\tilde g \ = \ z_{\tilde g} + {\frac {s_{\tilde g}}{n(n-1)}} \tilde g
$$
by using $s_g =s_{\tilde g}$ in the last step. 
 
Substituting the above expression for $ \tilde z_g$ and again using
$s_g = s_{\tilde g}$ in (\ref{sta}), we get
$$ { \frac {1}{(n-2)(n-3)}} z_{\tilde g} \owedge \tilde g=0.$$
This holds if and only if $z_{\tilde g}=0$, i.e.,  $\tilde g$ is Einstein.

Assuming that $\tilde g$ is Einstein, we now prove equality in (\ref{che}).
From the proof of Lemma \ref{nor}, it is clear that equality holds in
(\ref{che}) if and only $W_g(X,N,N,Y)=W_g(X,N,Z,Y)=0$ for any $X, \ Y, \ Z
\in T_p \Sigma$. Since $R_g(X,N,N,Y)=R_g(X,N,Z,Y)=0$ and
$R_g=S_g + W_g$, we just have to check that $S_g =0$ for these vectors. This
is a routine calculation involving the Kulkarni-Nomizu product which we skip.

Note that this remark proves that the product metric on $\C P^2 \times S^1$ 
satisfies $s \equiv 2 \sqrt 6 \vert W \vert $. \\

By combining Theorem
\ref{mai2} with the B\"ochner formula for harmonic 2-forms on a
4-manifold, we can now prove Corollary \ref{c1}.

\begin{proof}[Proof of Corollary \ref{c1}]
Suppose $dim H^2(\Sigma, \R) = b_2(\Sigma) \neq 0$. We will show that
$(\Sigma, h)$, where $h$ is the induced metric, is a self-dual K\"ahler surface
of non-negative scalar curvature.

We first show that $\Sigma$ cannot admit a metric
$h$ with $s_h >2 \sqrt 6 \vert W(h) \vert _h$. Suppose $h$ is such a metric.

By Hodge theory, $H^2( \Sigma, {\mathbb R}) \cong \H$, where $ \H$ is
the space of harmonic 2-forms. For four-manifolds we further have 
$\H =\H_{+} \oplus \H_{-}$, where $\H_{+}$ and $\H_{-}$ are the spaces
of self-dual and anti-self-dual harmonic 2-forms, respectively. 
Let $\omega \in \H_+$. The B\"ochner formula applied to $\omega$
(see, for instance, ~\cite{gur1}) gives
  $${\frac{1}{2}} \triangle \vert \omega \vert ^2 = \vert \nabla
\omega \vert ^2 -2W^+(\omega, \omega) +{\frac {1}{3}} s \vert \omega
\vert^2 ,$$ where $W=W^++W^-$ is the orthogonal decomposition of $W$
into self-dual and anti-self-dual parts.

 Since $W^+$ is a trace-free symmetric operator of the
three-dimensional vector bundle $ \Lambda_+$ (the bundle of self-dual
2-forms), we have, as in ~\cite{gur1}, $$-2W^+(\omega , \omega) \ge
-{\frac {2 \sqrt 6}{3}} \vert W^+ \vert \vert \omega \vert ^2. $$
Hence we get
\begin{equation}\label{boc}
{\frac{1}{2}} \triangle \vert \omega \vert ^2  \ \ge \ 
 \vert \nabla \omega \vert ^2 + {\frac {1}{3}} \bigl( s -2
\sqrt 6 \vert W^+ \vert  \bigr) \vert \omega \vert ^2  
  \  \ge  \ \vert \nabla \omega \vert ^2 + {\frac {1}{3}} \bigl( s -2
\sqrt 6 \vert W \vert  \bigr) \vert \omega \vert ^2  
\end{equation}

We have a similar inequality for any harmonic anti-self-dual 2-form.

Since $b_2 \neq 0$ by assumption, there is a non-zero harmonic 2-form
$\omega$, which we assume is self-dual without loss of generality. By
integrating the above inequality on $\Sigma$ we get a contradiction,
since we asumed that $s > 2 \sqrt 6 \vert W \vert $. 

Hence, if $b_2 \neq 0$, by Theorem \ref{mai2}, the induced metric $h$
on $\Sigma$ satisfies $s \equiv 2 \sqrt 6 \vert W \vert $. If $
\omega$ is as above, then, by (\ref{boc}), $ \nabla \omega =0$ and $s
\equiv 2 \sqrt 6 \vert W^+ \vert $. Hence the metric $h$ is K\"ahler and
$W^-=0$. If $\omega$ is
anti-self-dual, then we get $\nabla \omega =0$ and
$ W^+=0$. Hence, in either case, by reversing the orientation on
$\Sigma$ if necessary, $(\Sigma, h)$ is a
compact self-dual K\"ahler surface of non-negative scalar curvature.
\end{proof}

\end{document}